\documentclass{pjm}
\usepackage{oldgerm}
\usepackage{latexsym}
\usepackage{amsmath}
\usepackage{amssymb}
\usepackage{amsgen}
\usepackage{amscd} 
\usepackage{color}
\usepackage{epsfig}

\newtheorem{thm}{Theorem}[section]
\newtheorem{prop}[thm]{Proposition}
\newtheorem{cor}[thm]{Corollary}
\newtheorem{lem}[thm]{Lemma}

\theoremstyle{definition}
\newtheorem{dfn}[thm]{Definition}
\newtheorem{example}[thm]{Example}
\theoremstyle{remark}
\newtheorem{rem}[thm]{Remark}
\theoremstyle{remark}
\newtheorem{rmks}[thm]{Remarks}

\newcommand{\opr}[1]{\operatorname{#1}}
\newcommand{\transfer}{\opr{\stackrel{\bullet}{\phantom{.}}}}

\newcommand{\MDLM}{V(\mathbb{R})}

\newcommand{\Sharpo}[1]{\operatorname{\mathbf{S}}#1}

\newcommand{\phull}[1]{\opr{\diamond}(#1)} 
\newcommand{\lhull}[1]{\opr{\diamond}#1 }

\title{Quasidiagonal $C^\ast$-algebras and Nonstandard
Analysis}

\author{F. Javier Thayer} \address{9112 Decatur Ave S. \\
Bloomington, MN 55438} \email{jt@mitre.org}

\makeindex 

\begin{document} 
\begin{abstract}
Suppose $\mathbf{B}$ is an ultraproduct of finite dimensional
$C^\ast$-algebras.  We consider mapping and injectability properties
for separable $C^\ast$-algebras into $\mathbf{B}$. In the case of approximately
finite $C^\ast$-algebras, we obtain a classification of these mappings
up to inner conjugacy.  Using a Theorem of Voiculescu, we show that
for nuclear $C^\ast$-algebras injectability into an ultraproduct of
finite dimensional $C^\ast$-algebras is equivalent to quasidiagonality.
\end{abstract}

\maketitle
\section{Introduction}

In this paper we use nonstandard analysis~(\cite{albeverio-etal},
\cite{hurd-loeb}, \cite{keisler}) to investigate injectability into
$C^\ast$-algebras $\mathbf{B}$ which are infinitesimal hulls of
hyperfinite dimensional internal $C^\ast$-algebras $\mathfrak{B}$.
$\mathbf{B}$ is obtained from $\mathfrak{B}$ by considering the
subspace $\opr{Fin}(\mathfrak{B})$ of elements with norm $\ll \infty$
and identifying $x,y \in \opr{Fin}(\mathfrak{B})$ whenever $x - y$ has
infinitesimal norm.  Though $\mathbf{B}$ is a legitimate standard
$C^\ast$-algebra, it is very large, except in the uninteresting case
the original $\mathfrak{B}$ is finite dimensional.  We point out that
the $C^\ast$-algebras $\mathbf{B}$ are exactly ultraproducts of
finite-dimensional $C^\ast$-algebras (see~\cite{hensvino} or the
appendix).

A more interesting question from an operator theorist's viewpoint,
is which kinds of separable $C^\ast$-algebras are injectable into
$\mathbf{B}$ and what kinds of mappings exist from separable
$C^\ast$-algebras into $\mathbf{B}$.

We show the following: If $\mathbf{A}$ is an AF algebra,
Proposition~\ref{bijection-proposition} determines the inner conjugacy
classes of $C^\ast$-morphisms for $\mathbf{A}$ into a fixed
$\mathbf{B}$ in terms of certain projective systems of matrices with
nonnegative integer entries. Proposition~\ref{imbedding-condition}
gives a necessary and sufficient condition for injectability of an AF
algebra into a fixed $\mathbf{B}$ in similar terms.  Our approach uses
only the mapping properties of matrix algebras and a few basic
principles of nonstandard analysis.

In \S\ref{standard-characterization}, we give a characterization of
nuclear subhyperfinite $C^\ast$-algebras, that is nuclear
$C^\ast$-algebras which are injectable into some $\mathbf{B}$.
Theorem~\ref{characteriztion-of-qd} states that for nuclear
$C^\ast$-algebras, subhyperfiniteness is equivalent to
quasidiagonality.  To prove this characterization, we use a lifting
theorem for nuclear completely positive contractions reminiscent of
that of Choi and Effros~\cite{choi-effros}.  As a result, we obtain
the following characterization: A separable nuclear $C^\ast$-algebra
is quasidiagonal iff it is injectable into a countable ultraproduct of
finite dimensional $C^\ast$-algebras.  This characterization of
quasidiagonality is similar to that given by Theorem 5.2.2
of~\cite{blackadar-kirchberg} which considers injections into a
somewhat different kind of reduced product of finite dimensional
$C^\ast$-algebras.

There are numerous open questions we have not dealt with at this
time. Of particular interest is the precise relationship between
subhyperfinitess and nuclearity.  Another question, which is possibly
more of a set-theoretic nature, is whether the (enormous)
$C^\ast$-algebras $\mathbf{B}$ corresponding to distinct hyperfinite
dimensional internal $C^\ast$-algebras are non-isomorphic.

The structure of the paper is as follows:
\S\ref{preliminaries-section} is devoted to a review of basic ideas of
nonstandard analysis. In \S\ref{normed-algebra-section} we review list
general facts about internal $C^\ast$-algebras.  In
\S\ref{lifting-section}, we prove a number of basic lifting theorems
for projections and partial isometries in a somewhat unfamiliar
context, because the lifting property relates to a map between objects
in two different models of a theory.  However, the techniques used in
dealing with projections and partial isometries are much the same as
the techniques used by operator theorists in other situations.  These
lifting theorems are used in combination with various mapping
properties of hyperfinite dimensional algebras to obtain the basic
mapping properties for finite dimensional and AF algebras obtained in
\S\ref{matrix-sequences-section}.  Finally in
\S\ref{standard-characterization} we show the equivalence of
subhyperfiniteness and quasidiagonality in the nuclear case.  The
techniques used in this section are largely independent of those used
to investigate AF injectability.

The relationship between $C^\ast$-algebra theory and nonstandard
matrix algebras has been noticed before, notably
in~\cite{hinokuma-ozawa}.  The authors in this paper observe that the
infinitesimal hull of a matrix algebra is a $C^\ast$-algebra and
actually construct another related quotient space which turns out to
be a $\mathrm{II}_1$ factor.  However, there seems to be no suggestion
in that paper that $C^\ast$-algebras obtained from nonstandard matrix
algebras might have a useful relation to separable
$C^\ast$-algebras. An extensive and highly recommended reference for
quasidiagonal $C^\ast$-algebras is the paper~\cite{nate-brown}, which
records much recent progress in theory of quasidiagonal algebras
relevant to the present paper.  However the relation between
quasidiagonal algebras and ultraproducts of finite dimensional
dimensional algebras seems to be new.  In fact, the entire question of
injectability into such ultraproducts seems to be new as well.

\section{Preliminaries on Nonstandard
Analysis}\label{preliminaries-section}

We use a small amount of nonstandard analysis for which the first
pages of~\cite{keisler} suffice.  Given a set $\mathbf{G}$, the {\em
superstructure} over $\mathbf{G}$ is the set $V(\mathbf{G})$ defined
by:
$V_0(\mathbf{G}) = \mathbf{G}$, $V_{n+1}(\mathbf{G})
=V_{n}(\mathbf{G}) \cup \mathcal{P}(V_{n}(\mathbf{G}))$ and
$V(\mathbf{G}) = \bigcup_n V_{n}(\mathbf{G})$.
The main constituent of our working view of nonstandard analysis is a
map $\transfer: \MDLM \longrightarrow V(\transfer \mathbb{R})$ which
is the identity on $\mathbb{R}$ and satisfies the transfer principle;
c.f.~\cite{keisler} for details.  For the reader unfamiliar with
nonstandard analysis, we point out that a mapping satisfying the
transfer property can be shown to exist by using bounded
ultraproducts.  This is briefly reviewed in an appendix to this paper.
However, we recommend the discussion of the so-called ``extended
universe'' in~\cite{albeverio-etal} \S1.2.  Also note that the
transfer principle does not uniquely characterize the map $\transfer$.

\begin{rmks} We point out a notational difference between this paper
and most other papers in nonstandard analysis (including our
own~\cite{jt1}): the use of the prefixed symbol $\transfer$ instead of
the prefixed $\opr{\stackrel{\star}{\phantom{.}}}$ to denote the
transfer map.  This notational modification should prevent clashes
between the symbol used for the transfer map and the symbol used for
the adjoint mapping in involutive algebras.

For a similar reason we use the phrase ``involutive morphism'' in
place of the more frequently used expression ``$\star$-morphism''.
\end{rmks}

We use the modifier {\em standard} roughly, in referring to structures
in the universe of standard sets $\MDLM$.  The modifier {\em internal}
on the other hand refers to structures in the universe of internal
sets.  These are certain sets in $V(\transfer \mathbb{R})$
(see~\cite{keisler}, \S1 for a summary).  Strictly speaking, the
modifiers standard and internal should be used in reference to
particular models of a theory, but in practice the superstructure
approach to nonstandard analysis is sufficiently well-established to
allow us to gloss over these details. Thus we use without comment,
expressions such as standard or internal metric spaces, internal
groups etc.  We will also use the modifier {\em external} to draw
attention to the fact that a particular object is not internal or that
an assertion is used in an external context, for instance, {\em
external induction}.

We caution the experienced nonstandard analyst that we only use
countable saturation.

We use countable saturation as follows: If $\{x_n\}_{n \in
\mathbb{N}}$ is a sequence of internal elements of an internal set
$\mathcal{Y}$, then there is an internal family $\{\xi_\ell\}_{\ell
\in \transfer{\mathbb{N}}}$ which extends $\{x_n\}$.  

\begin{rmks} \label{remarks-on-saturartion} 
Saturation arguments are often used in conjunction with {\em
overspill} (also called {\em overflow}; see Proposition 1.2.7
of~\cite{albeverio-etal}) which states that any internal set which
contains $\mathbb{N}$ contains a hyperinteger interval $\{1,2, \ldots,
N\}$.  In particular, suppose $\{x_i\}_{i \in \mathbb{N}}$ is a
sequence in an internal set $A_1$ each term of which satisfies the
internal elementary statement $\mathbf{P}(x, A_1, \ldots , A_n)$,
where $A_i$ are parameters instantiated by fixed internal sets (see
\S1 of~\cite{keisler} for the definition of elementary
statement). Notice each term of the sequence is internal but the
sequence $\{x_i\}_{i \in \mathbb{N}}$ is external.  By countable
saturation, this sequence has an internal extension
$\{\xi_\ell\}_{\ell \in \transfer{\mathbb{N}}}$.  The set of $\ell \in
A_1$ for which $\mathbf{P}(\xi_\ell, A_1, \ldots , A_n)$ holds is
internal (requires the Internal Definition Principle) so by overspill,
there is an unlimited $N \in \transfer{\mathbb{N}}$ such that
$\mathbf{P}(\xi_\ell, A_1, \ldots , A_n)$ holds for $\ell \leq N$.
\end{rmks}

\subsection{Internal Normed
Spaces}\label{internal-normed-spaces-subsection}
See~\cite{albeverio-etal} for details on nonstandard analysis of
normed spaces.  If $\mathcal{E}$ is a $\transfer{\mathbb{C}}$ internal
normed space, then $\opr{Fin}(\mathcal{E})$ is the external set
$\{\phi \in \mathcal{E} : \|\phi\|_\mathcal{E} \ll \infty\}$, in other
words $\opr{Fin}(\mathcal{E})$ is the limited component of $0$. Note
that $\opr{Fin}(\mathcal{E})$ is a vector space over the field
$\mathbb{C}$. The mapping $\phi \mapsto \opr{st} \| \phi \|$ is a
($\mathbb{R}$-valued) seminorm on $\opr{Fin}(\mathcal{E})$.  This
seminorm factors through a norm on the quotient vector space
$\opr{Fin}(\mathcal{E})/\{\phi: \opr{st}\|\phi\| = 0\}$. The quotient
space with this norm is the {\em infinitesimal hull} of $\mathcal{E}$
and we denote it by $\Sharpo{(\mathcal{E})}$.  
%
%
We also denote the canonical quotient map $\opr{Fin}(\mathcal{E})
\rightarrow \Sharpo{(\mathcal{E})}$ by $\opr{\pi_\mathcal{E}}$ (or
$\opr{\pi}$ when the internal normed space $\mathcal{E}$ is evident
from context).  The norm on $\Sharpo{(\mathcal{E})}$ is characterized
by the property
$\|\opr{\pi_\mathcal{E}}(\phi)\| = \opr{st} \|\phi\|$.
It is well-known that if $\mathcal{E}$ is an internal normed space,
$\Sharpo{(\mathcal{E})}$ is a Banach space.  Note that if
$\mathcal{E}$ is hyperfinite dimensional, then $\mathcal{E}$ is
$\transfer{}$-complete. This means that any internal Cauchy sequence
$\{\phi_j\}_{j \in \transfer{\mathbb{N}}}$ is
$\transfer{}$-convergent.
\begin{rmks} \label{imbedding-in-collpase-transfer}
If $E$ is a standard normed space, there is a canonical linear
isometric map $E \rightarrow \Sharpo{(\transfer{E})}$. We will regard
this map is an inclusion.  This map is bijective iff $E$ is finite
dimensional. 

In particular, if $E$ is finite dimensional and $\phi \in
\transfer{E}$ is such that $\|\phi\| \ll \infty$ then
$\pi_{\mathcal{\transfer{E}}} (\phi)$ is an element of $E$.
$\pi(\phi)$ is usually referred to as the standard part of $\phi$.
Moreover, $\transfer{[\pi(\phi)]} \cong \phi$.
\end{rmks}
We adopt the following typographical convention to distinguish names
for internal normed spaces from those for standard normed spaces:
Internal spaces are denoted by calligraphic letters $\mathcal{E},
\mathcal{F}, \mathcal{H}$ while standard spaces are denoted by the
standard math fonts $E, F, H$.

We rely on similar typographical conventions to distinguish between
spaces of linear operators.  Suppose $\mathcal{E}, \mathcal{F}$ are
internal normed spaces.
$\mathfrak{L}(\mathcal{E},\mathcal{F})$ denotes the space of internal
$\transfer$-continuous linear maps $\mathcal{E} \rightarrow
\mathcal{F}$.  On the other hand, if $E$, $F$ are standard normed
spaces, $L(E,F)$ denotes the space of bounded linear maps $E
\rightarrow F$.
By transfer, there is an internal norm function $\| \cdot \|:
\mathfrak{L}(\mathcal{E},\mathcal{F}) \rightarrow
\transfer{\mathbb{R}}$.
$\mathfrak{L}_{\opr{Fin}}(\mathcal{E},\mathcal{F})$ consists of $T \in
\mathfrak{L}(\mathcal{E},\mathcal{F})$ such that the internal norm
$\|T\|$ is limited.  If $T \in
\mathfrak{L}_{\opr{Fin}}(\mathcal{E},\mathcal{F})$ then $T$ maps
$\opr{Fin}(\mathcal{E})$ into $\opr{Fin}(\mathcal{F})$ and $T \,|
\,\opr{Fin}(\mathcal{E})$ factors through a bounded linear map
$\Sharpo{(T)}: \Sharpo{(\mathcal{E})} \rightarrow
\Sharpo{(\mathcal{F})}$.  $\Sharpo{(T)}$ is characterized by
\begin{equation}\label{sharpo-characterization}
\Sharpo{(T)}(\opr{\pi}{(\phi)}) = \opr{\pi}{(T \phi)} \quad \mbox{for $\phi \in
\mathcal{E}$.}
\end{equation}
If $T \in \mathfrak{L}(\mathcal{E},\mathcal{F})$ is isometric, then
$\Sharpo{(T)}$ is isometric.

$\mathcal{E} \mapsto \Sharpo{(\mathcal{E})}$ is a functor
on the category of internal normed spaces and maps of limited norm
into the category of Banach spaces and bounded linear maps.  

\begin{rmks}
If $T \in {\mathfrak{L}_{\opr{Fin}}(\mathcal{E},\mathcal{F})}$, the
operator $\Sharpo{(T)}: \Sharpo{(\mathcal{E})} \rightarrow
\Sharpo{(\mathcal{F})}$ can also be interpreted as the image of $T$ in
$\Sharpo{(\mathfrak{L}(\mathcal{E},\mathcal{F}))}$ under the quotient
map $\opr{\pi}:\mathfrak{L}_{\opr{Fin}}(\mathcal{E},\mathcal{F})
\rightarrow \Sharpo{(\mathfrak{L}(\mathcal{E},\mathcal{F}))}$.
More precisely, 
the map $\Sharpo{}:
\mathfrak{L}_{\opr{Fin}}(\mathcal{E},\mathcal{\mathcal{F}})
\rightarrow L(\Sharpo{(\mathcal{E})}, \Sharpo{(\mathcal{F})})$ factors
through a linear isometric map $\Phi:
\Sharpo{(\mathfrak{L}(\mathcal{E},\mathcal{F}))} \rightarrow
L(\Sharpo{(\mathcal{E})}, \Sharpo{(\mathcal{F})})$.
The isometric property of $\Phi$ means $\opr{st} \| T \| =
\|\Sharpo{(T)}\|$ for $T \in
\mathfrak{L}_{\opr{Fin}}(\mathcal{E},\mathcal{\mathcal{F}})$. To
verify this, note that by~(\ref{sharpo-characterization}),
$$\|\Sharpo{(T)} \opr{\pi} \phi \| = \opr{st} \| T \phi\| \leq \opr{st}
\| T \| \ \|\opr{\pi} \phi\|,$$
and thus $ \|\Sharpo{(T)}\| \leq \opr{st} \| T \|$.  On the other
hand, let $\epsilon \cong 0$ be a positive hyperreal, $\phi \in
\mathcal{E}$ such that $\|\phi\| = 1$ and $\| T \phi \| \geq \| T \| -
\epsilon$.
Then
$$\opr{st} \| T \| = \opr{st} (\| T \|- \epsilon) \leq \opr{st} \| T
\phi \| = \|\Sharpo{(T)} \opr{\pi}(\phi)\| \leq \|\Sharpo{(T)}
\|.$$
\end{rmks}

\begin{rmks} $\opr{\pi}$ maps part of an internal normed space
$\mathcal{E}$ linearly into the Banach space $\Sharpo{(\mathcal{E})}$.
Linearity of $\opr{\pi}$ means
$\opr{\pi}{\bigl(\sum_i x_i\bigr)} = \sum_i \opr{\pi}{(x_i)}$
for {\em finite} families $\{x_i\}$. Observe there is no reasonable
way to interpret the RHS of this formula if the index set is
hyperfinite but not finite.  A similar comment applies to the map $T
\mapsto \Sharpo{(T)}$.
\end{rmks}

\begin{example}\label{lp-spaces-example}
Suppose $\mathfrak{X}$ is hyperfinite and $1 \leq p < \infty$ in
$\transfer{\mathbb{R}}$. $\mathcal{L}^p(\mathfrak{X})$ is the internal
vector space of internal functions $\psi: \mathfrak{X} \rightarrow
\transfer{\mathbb{C}}$ equipped with the norm
$\|\psi\|_p=\sqrt[p]{\sum_{x \in \mathfrak{X}} |\psi(x)|^p}$.
\end{example}

If $\mathcal{H}$ is an inner product space, then by transfer, the
Cauchy-Schwartz inequality is valid in $\mathcal{H}$.  From this
follows that the inner product on $\mathcal{H}$ factors through an
inner product on the Banach space $\Sharpo{(\mathcal{H})}$.  In
particular $\mathcal{H} \mapsto \Sharpo{(\mathcal{H})}$ is also a
functor on the category of internal inner product spaces and maps of
limited norm into the category of Hilbert spaces and bounded linear
maps.

Note that if $\mathcal{H}$, $\mathcal{K}$ are hyperfinite dimensional
inner product spaces and $T:\mathcal{H} \rightarrow \mathcal{K}$ is a
linear map then the adjoint $T^\ast$ is always defined and has the
property that $\|T^\ast\| =\|T\|$.  If $T$ has limited norm,
$\Sharpo{(T^\ast)} = {\Sharpo{(T)}}^\ast$.

\begin{example}
If $\mathfrak{X}$ is hyperfinite, \label{l2-spaces-example}
$\mathcal{L}^2(\mathfrak{X})$ has an internal inner product
$$\langle \psi , \phi \rangle = \sum_{x \in \mathfrak{X}} \psi(x)
\overline{\phi(x)},$$
and clearly $\|\psi\|_2 = \sqrt{\langle \psi , \psi \rangle}$.
We alert the reader to the fact that unless the internal
dimension of $\mathcal{H}$ an element of $\mathbb{N}$,
$\Sharpo{(\mathcal{H})}$ is highly nonseparable. \end{example}

If $K$ is a separable Hilbert space, $\transfer{K}$ has an internal
orthonormal basis indexed on $\transfer{\mathbb{N}}$.

\section{Internal $C^\ast$-Algebras}  \label{normed-algebra-section}

Our typographical strategy for distinguishing internal from standard
$C^\ast$-algebras is as follows: boldface fonts $\mathbf{A}$,
$\mathbf{B}$ denote standard $C^\ast$-algebras, fraktur fonts
$\mathfrak{A}$, $\mathfrak{B}$, $\mathfrak{M}$ denote internal ones.

\begin{rem} \label{algebras-are-unital}
In this paper, we only consider algebras (whether internal or
external) with unit.  All involutive morphisms between algebras (in
particular, all inclusions) will be assumed to be unit preserving.
\end{rem}

Suppose $\mathfrak{A}$ is an internal normed algebra over the field
$\transfer{\mathbb{C}}$.  Note that $\mathfrak{A}$ is an algebra over
the field $\transfer{\mathbb{C}}$ and $\opr{Fin}(\mathfrak{A})$ is an
algebra over $\mathbb{C}$.  The Banach space $\Sharpo{(\mathfrak{A})}$
can be turned into an algebra such that $\pi:\opr{Fin}(\mathfrak{A})
\rightarrow \Sharpo{(\mathfrak{A})}$ is a ring homomorphism.  This
follows from the fact $\{x \in \mathfrak{A}: \|x\| \cong 0\}$ is a two
sided ideal in $\opr{Fin}(\mathfrak{A})$.  If $\mathfrak{A}$ is an
internal involutive normed algebra with the property $\|x \, x^\ast \|
= \|x\|^2$, then $\Sharpo{(\mathfrak{A})}$ is a $C^\ast$-algebra:
$$\|\opr{\pi}{(x)} \, \opr{\pi}{(x)}^\ast \| = \|\opr{\pi}{(x \, x^\ast)} \| =
\opr{st} \|x \, x^\ast\| = \opr{st} \|x\|^2 = \|\opr{\pi}{(x)}\|^2.$$

The facts listed in the following examples are well-known and are
stated without proof.

\begin{example} \label{algebra-of-operators}
If $\mathcal{H}$ is an internal inner product space,
$\mathfrak{L}(\mathcal{H})$ is an internal $C^\ast$-algebra.  This
follows by transfer.  $\mathfrak{L}_{\opr{Fin}}(\mathcal{H})$ is an
(external) algebra over $\mathbb{C}$ closed under adjoints.  The map
$T \mapsto \Sharpo{(T)}$ is an involutive algebra homomorphism
$\mathfrak{L}_{\opr{Fin}}(\mathcal{H}) \rightarrow
L(\Sharpo{(\mathcal{H})})$.  The image of this map is a
$C^\ast$-algebra.
\end{example}

\begin{example}
If $\mathbf{A}$ is a standard $C^\ast$-algebra,
$\transfer{\mathbf{A}}$ is an internal $C^\ast$-algebra.  The
embedding $\mathbf{A} \rightarrow \Sharpo{(\transfer{\mathbf{A}})}$ is
an injective morphism of $C^\ast$-algebras.  This map is surjective
iff $\mathbf{A}$ is finite dimensional.

Note that any separable $C^\ast$-algebra is isomorphic to a
$C^\ast$-algebra in the universe $\MDLM$.  $\MDLM$ of course also
contains enormously large $C^\ast$-algebras such as countable
ultrapowers of separable ones.  It is easy to show however, $\MDLM$
does not contain copies of every $C^\ast$-algebra.
\end{example}

\begin{example}
If $\mathcal{H}$ is hyperfinite dimensional, the internal
$C^\ast$-algebra $\mathfrak{L}(\mathcal{H})$ can also be viewed as an
algebra of square {\em internal} $n \times n$ matrices with entries in
$\transfer{\mathbb{C}}$.
The size $n$ of the matrix algebra is the internal dimension of
$\mathcal{H}$.
Note that the elements of $\opr{M}_n(\transfer{\mathbb{C}})$ are {\em
internal} functions
$$\{1, 2, \ldots , n\} \times \{1, 2, \ldots , n\} \rightarrow
\transfer{\mathbb{C}}.$$
\end{example}

\begin{example}
Suppose $\{\mathfrak{A}_i\}_{i \in F}$ is a hyperfinite
internal family of $C^\ast$-algebras.  Then
$\bigoplus_{i \in F} \mathfrak{A}_i$ is the space of internal families
$\bar{x}= \{x_i\}_{i \in F}$ such that $x_i \in \mathfrak{A}_i$,
equipped with the pointwise algebraic operations and the norm
$\|\bar{x}\| = \max_i \|x_i\|_{\mathfrak{A}_i}$. $\bigoplus_{i \in F}
\mathfrak{A}_i$ is an internal $C^\ast$-algebra.
\end{example}

\begin{example}
If $n \in \transfer{\mathbb{N}}$ and $\mathfrak{A}$ is an internal
$C^\ast$-algebra, $\opr{M}_n(\mathfrak{A})$, the algebra of internal
$n \times n$ with entries in $\mathfrak{A}$ is a $C^\ast$-algebra.
Note that if $n \in \mathbb{N}$, then
\begin{equation}\label{M-Sharpo-is-Sharpo-M}
\Sharpo{(\opr{M}_n(\mathfrak{A}))} \cong
\opr{M}_n\bigl(\Sharpo{(\mathfrak{A})}\bigr).
\end{equation}
More precisely, the map 
$\opr{M}_n(\pi): \opr{M}_n(\mathfrak{A}) \rightarrow
\opr{M}_n\bigl(\Sharpo{(\mathfrak{A})\bigr)}$ is naturally equivalent to
the infinitesimal identification map $\pi_{\opr{M}_n(\mathfrak{A})}$,
that is there is a commutative diagram
\begin{equation}\label{matricial-equivalence-diagram}
\begin{CD}
\opr{M_n}(\mathfrak{A}) @>{{\opr{id}_{\opr{M_n}(\mathfrak{A})}}}>>
\opr{M_n}(\mathfrak{A}) \\
@V{\opr{M_n}}(\pi)VV @VV\pi_{\opr{M_n}(\mathfrak{A})}V \\
\opr{M_n}(\Sharpo{(\mathfrak{A})}) @>{\cong}>>
\Sharpo{(\opr{M_n}(\mathfrak{A}))}
\end{CD}
\end{equation}
In the case $n$ is a nonstandard hyperinteger, the RHS
of~(\ref{M-Sharpo-is-Sharpo-M}) is meaningless. \end{example}

\subsection{Hyperfinite Dimensional
$C^\ast$-algebras}\label{fd-category}
The structure of finite dimensional $C^\ast$-algebras and involutive
morphisms between them is well-known (See~\cite{jones-sunder}, \S3.2
or~\cite{davidson}, Chapter III). The corresponding results for
hyperfinite dimensional $C^\ast$-algebras, which we now state, follow
by transfer.

A hyperfinite dimensional $C^\ast$-algebra $\mathfrak{A}$ is
canonically isomorphic to a hyperfinite direct sum $\mathfrak{A} =
\bigoplus_{e \in \opr{min}\mathfrak{A} } \mathfrak{A} e$, where
$\opr{min}\mathfrak{A}$ is the set of minimal nonzero central
projections of $\mathfrak{A}$. Each $C^\ast$-algebra $\mathfrak{A} e$
is internally isomorphic (though in a noncanonical way) to the full
matrix algebra $\opr{M}_{\opr{dim}(e)}(\transfer{\mathbb{C}})$. The
hyperfinite family indexed on $\opr{min}\mathfrak{A}$ given by
$\overline{\opr{dim}}(\mathfrak{A})= \{\opr{dim}(e)\}_{e \in
\opr{min}\mathfrak{A}}$ is called the {\em dimension vector} of
$\mathfrak{A}$.

If $\mathfrak{A}$, $\mathfrak{B}$ are hyperfinite
dimensional $C^\ast$-algebras, to any involutive morphism
$h:\mathfrak{A} \rightarrow \mathfrak{B}$
corresponds an internal family $\Lambda$ of nonnegative hyperintegers
indexed on $\opr{min} \mathfrak{B} \times \opr{min} \mathfrak{A}$.  This
family is defined as follows:
If $(f,e) \in \opr{min} \mathfrak{B} \times \opr{min} \mathfrak{A}$,
then $\Lambda_{f e}$ is the multiplicity of the irreducible
representation of $\mathfrak{A}$ corresponding to the projection $e$
in the representation $x \mapsto h(x) f$.  The mapping $h \mapsto
\Lambda(h)$ is internal.  We will refer to $\Lambda(h)$ as the {\em
mapping matrix} of $h$ (or {\em inclusion matrix} in case $h$ is
injective).

The following properties are well-known:
\begin{enumerate}
\item \label{dimension-matrix-equation} The following matrix equation
holds:
$\overline{\opr{dim}}(\mathfrak{B}) = \Lambda \
\overline{\opr{dim}}(\mathfrak{A})$.
\vspace{1mm}
\item If $h:\mathfrak{A} \rightarrow \mathfrak{B}$ and $g:
\mathfrak{B} \rightarrow \mathfrak{C}$ are involutive morphisms, then
$\Lambda(g \circ h) = \Lambda(g) \Lambda(h)$.
\vspace{1mm}
\item \label{same-mapping-matrix-equivalent} $h:\mathfrak{A}
\rightarrow \mathfrak{B}$, $h': \mathfrak{A} \rightarrow
\mathfrak{B}'$ have the same mapping matrix iff there is an internal
$C^\ast$-algebra isomorphism $\Phi: \mathfrak{B}\rightarrow
\mathfrak{B}'$ such that $h' = \Phi \circ h$.
\vspace{1mm}
\item \label{unital-morphism-item} For each $f \in \opr{min} \mathfrak{B}$
there is an $e \in \opr{min} \mathfrak{A}$ such that $\Lambda_{f e} >
0$. This follows immediately from the fact $h$ is unital.
\vspace{1mm}
\item If $\mathfrak{A}$, $\mathfrak{B}$ are hyperfinite dimensional
$C^\ast$-algebras and $\Lambda$ is a family of nonnegative
hyperintegers indexed on $\opr{min} \mathfrak{B} \times \opr{min}
\mathfrak{A}$ satisfying the condition in item~\ref{unital-morphism-item},
then there is an involutive morphism $h:\mathfrak{A} \rightarrow
\mathfrak{B}$ for which $\Lambda$ is the mapping matrix.
\vspace{1mm}
\item \label{injective-mapping-matrix-property} $h$ is injective iff
for each $e \in \opr{min} \mathfrak{A}$ there is an $f \in \opr{min}
\mathfrak{B}$ such that $\Lambda_{f e} > 0$.
\item \label{dimension-matrix-equation-converse} Conversely, if
$\Lambda$ is a  matrix such that items
~\ref{dimension-matrix-equation}) and ~\ref{unital-morphism-item}) hold, 
then there is a unital morphism $h:\mathfrak{A} \rightarrow
\mathfrak{B}$ such that $\Lambda(h) = \Lambda$.
\end{enumerate}

%
%
\begin{dfn}
A family of nonnegative hyperintegers $\{\Lambda_{j \ i} \}_{(j \ i)
\in J \times I}$ such that for each $f \in J$ there is an $e \in I$
such that $\Lambda_{f e} > 0$ is called a mapping matrix.  If for
each $e \in I$ there is an $f \in J$ such that $\Lambda_{f e} > 0$,
$\Lambda$ is called an inclusion matrix.
\end{dfn}
\subsection{Functional Calculus}
The existence of a $\transfer$-continuous functional calculus in an
internal $C^\ast$-algebra for selfadjoint elements follows from the
existence of the continuous functional calculus for selfadjoint
elements in standard $C^\ast$-algebras. For example, any internal
function $f$ such that for any $\theta \in \opr{dom} f$, $f \opr{|}
[\theta - \delta, \theta] \cap \opr{dom}f $ and $f \opr{|} [\theta,
\theta + \delta]\cap \opr{dom}f$ are polynomial functions for some
$\delta > 0$ are $\transfer$-continuous. For another example, if $f$
is a standard continuous function, $\transfer{f}$ is
$\transfer$-continuous.

We now observe that the internal functional calculus on $\mathfrak{A}$
is compatible with the standard one of $\Sharpo{(\mathfrak{A})}$, as
stated in the following:

\begin{lem} \label{compatibility-of-functional-calc-prop}
Suppose $x \in \opr{Fin}(\mathfrak{A})$ is self-adjoint. If
$f:\mathbb{R} \rightarrow \mathbb{R}$ is continuous, then
\begin{equation}\label{compatibility-of-functional-calc}
\opr{\pi}{([\transfer{f}](x))} = f(\opr{\pi}{(x)}).
\end{equation}
\end{lem}
\begin{proof}
Without loss of generality we can assume $\|x\| \leq 1$.  Since $\pi$
is a homomorphsim, equation~(\ref{compatibility-of-functional-calc})
holds for standard polynomials $f$. The RHS
of~(\ref{compatibility-of-functional-calc}) is norm continuous in $f
\in \opr{C}[-1,1]$.  Moreover, by transfer, it follows that $f \mapsto
f(x)$ is an internal linear mapping of norm $\leq 1$ as $f$ ranges
through the space of $\transfer$-continuous functions on
$\transfer[-1,1]$. Thus
$f \mapsto \opr{\pi}{(\transfer{f}(x))}$ is a linear mapping of norm
$\leq 1$ as $f$ ranges through the space of continuous functions on
$[-1,1]$.  The result follows by continuity and Stone-Weierstrauss.
\end{proof}

\section{Lifting Projections and Partial Isometries}
\label{lifting-section} 
Let $\mathfrak{A}$ be an internal $C^\ast$-algebra.  The relation
between $\mathfrak{A}$ and $\Sharpo{(\mathfrak{A})}$ bears some formal
resemblance to that between the standard $C^\ast$-algebra $L(H)$ and
the Calkin Algebra $L(H)/K(H)$, particularly as regards projections
and partial isometries.
Projection in a $C^\ast$-algebra (whether standard or internal) means
self-adjoint projection and ``$\preceq$'' denotes the ordering
relation on projections.
\begin{prop}\label{lifting-partial-isometries-etc}
If $w\in\Sharpo{(\mathfrak{A})}$ is a partial isometry, then there is a
partial isometry $\tilde{w} \in \mathfrak{A}$ such that $\opr{\pi}{(\tilde{w})} = w$.
Moreover, if $p,p' \in \mathfrak{A}$ are projections such
that
$$\opr{initial \ projection}w \preceq \opr{\pi}{(p)} \mbox{
and } \opr{final \ projection }w \preceq \opr{\pi}{(p')},$$
$\tilde{w}$ can be taken so that its initial projection $\preceq p$ and its
final projection  $\preceq p'$.
\end{prop}
\begin{proof} 
Let $b \in \mathfrak{A}$ such that $\opr{\pi}{(b)} = w$. Replace $b$ by
$p' \ b \ p$ and let $x =b^\ast b$.  Thus
$\opr{\pi}{(x)} = w^\ast w = e$ and $x^2 - x \cong 0$.
Let $h:\mathbb{R} \rightarrow
\mathbb{R}$ be the function 
$$h(\theta) = \left\{\begin{array}{cc}
\theta^{-1/2} & \mbox{ if $\theta \geq 2/3$}, \\
0 & \mbox{ if $\theta \leq 1/3$},\end{array}\right.$$
and is linear in $[1/3, 2/3]$. By the internal functional calculus,
$[\transfer{h}](x)^2 \ x$ is a projection.  By
Lemma~\ref{compatibility-of-functional-calc-prop},
$\opr{\pi}{([\transfer{h}](x))} = h\bigl(\opr{\pi}{(x)}\bigr) = h(e) =
e$.
$\tilde{w}=b \,  [\transfer{h}](x)$ is a partial isometry:
$$(b \, [\transfer{h}](x))^\ast \, b \, [\transfer{h}](x) = [\transfer{h}](x) \, b^\ast b \, [\transfer{h}](x) =  {[\transfer{h}](x)}^2 \, x.$$
Since $b$ has initial projection $\preceq p$ and final projection is
$\preceq p'$, it follows that the same is true for $b \,
[\transfer{h}](x)$.  Finally,
$\opr{\pi}{(\tilde{w})} = \opr{\pi}{(b)}\,
\opr{\pi}{([\transfer{h}](x))} = w \, e = w$.\end{proof}

\begin{rem}\label{no-index-obstructions}
Unlike the situation with Calkin map $L(H) \rightarrow L(H)/K(H)$
there are no index obstructions to lifting mappings from finite
dimensional algebras into $\Sharpo{(\mathfrak{A})}$. The relevant fact
is that if $p \in \mathfrak{A}$ is a projection such that $\pi(p)=1$,
then $p=1$.  To see this, note that $1 - p$ is also a projection and
$1 - p \cong 0$.  By the Gelfand isomorphism the only projection with
norm $<1$ in a standard $C^\ast$-algebra is $0$.  By transfer, the
same fact is true for projections in internal $C^\ast$-algebras.  Thus
$1 - p = 0$.  In particular, by (external) induction on $r \in
\mathbb{N}$:

Suppose $\{p_k\}_{1 \leq k \leq r}$ (with $r \in \mathbb{N}$) are
orthogonal projections in the internal $C^\ast$-algebra $\mathfrak{A}$. If 
$\sum_{\ell=1}^r \pi(p_\ell) = 1$, then $\sum_{\ell=1}^r p_\ell =1$.
\end{rem}

A basic notion for dealing with finite-dimensional $C^\ast$-algebras
are {\em matrix units}. We follow the definitions and the notation of
\S7.1 of~\cite{rordam-larsen-laustsen}, with the following two
caveats: {\em all our matrix units are nonzero} and the matrix units
are not required to generate the $C^\ast$-algebra.
\begin{dfn}
A system of matrix units for $\mathbf{A}$ is a family $\{e_{i
j}^k\}_{1 \leq k \leq r, 1 \leq i, j \leq n_k}$ of nonzero elements of
$\mathbf{A}$ which satisifies (i), (ii) and (iii) of \S7.1
of~\cite{rordam-larsen-laustsen} and
$$1_{\mathbf{A}} = \sum_{k=1}^r \sum_{i=1}^{n_k} e_{i i}^k.$$
We also add the following bit of jargon: the family of positive
integers $\{n_k\}_{1 \leq k \leq r}$ is the dimension vector of
$\{e_{i j}^k\}$.
\end{dfn}  

We will freely use basic facts about matrix units, most notably the
correspondence between, on the one hand the set of injective
involutive morphisms $\mathbf{A} \rightarrow \mathbf{B}$ and on the
other the set of pairs consisting of matrix units $\{e_{i j}^k\}$
which span $\mathbf{A}$ (so that $\mathbf{A}$ is finite dimensional)
and matrix units $\{f_{i j}^k\}$ in $\mathbf{B}$ with the same
dimension vector as $\{e_{i j}^k\}$. These facts are well known, and
again we refer the reader to 7.1 of~\cite{rordam-larsen-laustsen} for
details.

\begin{rem}\label{reduced-matrix-units}
Note that to specify a system of matrix units with dimension vector
$\{n_k\}_{1 \leq k \leq r}$ for $\mathbf{A}$ it suffices to provide a
a {\em reduced system of matrix units}, that is a sequence of pairwise
orthogonal projections $\{p_k\}_{1 \leq k \leq r}$ such that
$1_{\mathbf{A}} = \sum_{k=1}^r p_k$ and family of partial isometries
$\{u_j^k\}_{1 \leq k \leq r, 1 \leq j \leq n_k}$ such that $1 \leq k \leq r$
$$
\begin{aligned}
\sum_{j=1}^{n_k} \opr{final \ projection} u_j^k & = p_k, \\
\opr{initial \ projection} u^k_j & = u^k_1 \quad 1 \leq j \leq n_k.
\end{aligned}
$$
\end{rem}

Clearly there is a corresponding internal notion of system of matrix
units and as always we will not hesitate to use the corresponding
transferred facts about matrix units.

As a corollary to the lifting theorem for partial isometries and the
above remark:

\begin{cor} 
Suppose $\mathfrak{A}$ is an internal $C^\ast$-algebra and $\{e_{i
j}^k\}$ a system of matrix units in $\Sharpo{(\mathfrak{A})}$.  Then
there is an system of matrix units $\{\phi_{i j}^k\}$ in
$\mathfrak{A}$ with the same dimension vector as $\{e_{i j}^k\}$ for
which
$\pi(\phi_{i j}^k) = e_{i j}^k$.
\end{cor}
%
%
\begin{rem} If $\mathbf{A}$ is a standard $C^\ast$-algebra and $\{e_{i
j}^k\}$ is a system of matrix units in $\mathbf{A}$, then
$\{\transfer{e_{i j}^k}\}$ is a system of matrix units in
$\transfer{\mathbf{A}}$.  If we regard $\mathbf{A}$ as isometrically
embedded in $\Sharpo{(\transfer{\mathbf{A}})}$ (see
Remark~\ref{imbedding-in-collpase-transfer}), then $\pi(\transfer{e_{i
j}^k}) = e_{i j}^k$. Thus in the special case
$\mathfrak{A}=\transfer{\mathbf{A}}$ lifting matrix units is
trivial.
\end{rem}

\begin{cor}\label{lifting-individual-morphisms}
Suppose $\mathbf{A}$ is a standard finite dimensional
$C^\ast$-algebra, $\mathfrak{B}$ an internal $C^\ast$-algebra,
$\phi:\mathbf{A} \rightarrow \Sharpo{(\mathfrak{B})}$ an involutive
morphism.  Then there is an internal involutive morphism
$\Phi:\transfer{\mathbf{A}} \rightarrow \mathfrak{B}$ of internal
$C^\ast$-algebras which makes the diagram
\begin{equation}\label{morphism-lifting-diagram}
\begin{CD}
\transfer{\mathbf{A}}  @>{\Phi}>>   \mathfrak{B} \\
@A{\transfer{}}AA     @VV{\pi}V \\
\mathbf{A}   @>{\phi}>>   \Sharpo{(\mathfrak{B})}
\end{CD}
\end{equation}
commutative. In particular, $\Sharpo{(\Phi)} = \phi$.
\end{cor}
\begin{proof}
Without loss of generality we can assume $\phi$ is injective.  Consider
a system of matrix units $\{e_{i j}^k\}_{1 \leq k \leq r, 1 \leq i, j
\leq n_k}$ for $\mathbf{A}$. Since $\phi$ is injective, $\{\phi(e_{i
j}^k)\}_{1 \leq k \leq r, 1 \leq i, j \leq n_k}$ is a system of matrix
units in $\mathfrak{B}$. Extend by $\mathbb{C}$-linearity.
\end{proof}
\begin{prop}\label{lifting-sequences-of-moorphisms}
Suppose $\{\mathbf{A}_i\}_{i \in \mathbb{N}}$ is a nondecreasing
sequence of finite dimensional $C^\ast$-subalgebras of $\mathbf{A}$
such that $\bigcup_{\ell=1}^\infty \mathbf{A}_\ell$ is norm dense in
$\mathbf{A}$. Suppose also $\mathfrak{B}$ is an internal
$C^\ast$-algebra and $\{\phi_\ell\}_{\ell \in \mathbb{N}}$ is a
sequence of involutive morphisms $\mathbf{A}_\ell \rightarrow
\Sharpo{(\mathfrak{B})}$ such that
$$\phi_{\ell+1} \opr{|} \mathbf{A}_\ell = \phi_\ell \quad \mbox{for $\ell \in
\mathbb{N}$.}$$
Then there is a hyperfinite dimensional $C^\ast$-algebra $\mathfrak{A}
\subseteq \transfer{\mathbf{A}}$ such that $\transfer{\mathbf{A}_\ell}
\subseteq \mathfrak{A}$ for all $\ell \in \mathbb{N}$ and an internal
involutive morphism $\Phi: \mathfrak{A} \rightarrow \mathfrak{B}$ such
that
\begin{equation}\label{restriction-equation}
\Sharpo{(\Phi)} \opr{|} \mathbf{A}_\ell = \phi_\ell\quad \mbox{for
$\ell \in \mathbb{N}$,}
\end{equation}
where we view $\mathbf{A}_\ell$ as isometrically imbedded in
$\Sharpo{(\mathfrak{A})}$ as is justified by
Remark~\ref{imbedding-in-collpase-transfer}.
\end{prop}
\begin{proof}
Consider the sequence of internal $C^\ast$-algebras
$\transfer{\mathbf{A}}_\ell = \mathfrak{A}_\ell$.  Observe that
$\mathfrak{A}_\ell \subseteq \mathfrak{A}_k \subseteq
\transfer{\mathbf{A}}$ for $k \geq \ell$.  By
Corollary~\ref{lifting-individual-morphisms}, and the fact each
$\mathbf{A}_\ell$ is finite-dimensional there is an internal
involutive morphism $\Phi_\ell: \mathfrak{A}_\ell \rightarrow
\mathfrak{A}$ such that
\begin{equation} \label{commutative-diagram-of-lifting}
\begin{CD}
\transfer{\mathbf{A}_\ell}  @>{\Phi_\ell}>>   \mathfrak{B} \\
@A{\transfer{}}AA     @VV{\pi}V \\
\mathbf{A}_\ell   @>{\phi_\ell}>>   \Sharpo{(\mathfrak{B})}
\end{CD}
\end{equation}	
commutes.
By saturation and overspill, the sequences $\{\Phi_\ell\}_{\ell \in
\mathbb{N}}$ and $\{\mathfrak{A}_\ell\}_{\ell \in \mathbb{N}}$ extend
to internal hyperfinite sequences $\{\Phi_\ell\}_{1 \leq \ell \leq
N}$, $\{\mathfrak{A}_\ell\}_{1 \leq \ell \leq N}$ such that for $\ell
\leq N$, $\Phi_\ell$ is an involutive morphism defined on the
hyperfinite dimensional $C^\ast$-algebra $\mathfrak{A}_\ell$.  If $x
\in \mathfrak{A}_\ell$ with $\|x\| \ll \infty$ and $\ell \leq k \in
\mathbb{N}$,
$$\Phi_\ell(x) \cong \Phi_\ell(\transfer{[\pi_{\mathfrak{A}_\ell}(x)]}) \cong
\Phi_k (\transfer{[\pi_{\mathfrak{A}_k}(x)]}) \cong
\Phi_k(x).$$
By the argument used to prove Robinson's lemma (see~\cite{nelsonprob},
Theorem 5.5), we conclude there is an unlimited $M$ such that
$\Phi_\ell(x) \cong \Phi_M(x)$ for $\ell \leq M$ and every $x \in
\mathfrak{A}_\ell$ of limited norm.  For completeness, we spell out
the details: for every $k \in \mathbb{N}$, $\|\Phi_\ell(x) -
\Phi_k(x)\| \leq 1/k$ for all $\ell \leq k$ and all $x \in
\mathfrak{A}_\ell$ with $\|x\| \leq 1$.  By overspill, there is an $M
\cong \infty$ such that $\|\Phi_\ell(x) - \Phi_M(x)\| \leq 1/M \cong
0$ for every $\ell \leq M$ and every $x \in \mathfrak{A}_\ell$ with
$\|x\| \leq 1$.

Therefore, for $\ell \leq M$, 
$$\Sharpo{(\Phi_M)} \opr{|} \mathbf{A}_\ell = \Sharpo{(\Phi_\ell)} =
\phi_\ell.$$\end{proof}
It immediately follows that if $\mathbf{A}$ is an approximately finite
$C^\ast$-algebra of $\Sharpo{(\mathfrak{B})}$, then there is a
hyperfinite dimensional $C^\ast$-algebra $\mathfrak{A} \subseteq
\mathfrak{B}$ such that $\pi (\mathfrak{A}) \supseteq \mathbf{A}$.
This however is a special case of
Proposition~\ref{characterization-of-subhyperfinite}.

\section{Matrix Sequences} \label{matrix-sequences-section}

Involutive morphisms $\phi, \psi: \mathbf{A} \rightarrow \mathbf{B}$
are {\em inner conjugate} iff there is a unitary $u \in \mathbf{B}$
intertwining $\phi$, $\psi$, i.e., for all $x \in \mathbf{A}$,
$\phi(x) = u^\ast \psi(x) u$.  We now prove the following peculiar
compactness property:

\begin{prop} \label{conjugacy-condition}
Suppose $\mathfrak{B}$ is an internal $C^\ast$-algebra, $\mathbf{A} =
\overline{\bigcup_{i=1}^\infty \mathbf{A}_i}$ is separable, where
$\{\mathbf{A}_i\}_{i \in \mathbb{N}}$ is a nondecreasing sequence of
$C^\ast$-algebras and $\phi, \psi:\mathbf{A} \rightarrow
\Sharpo{(\mathfrak{B})}$ are $C^\ast$-algebra morphisms.  If for all
$i \in \mathbb{N}$, $\phi \opr{|} \mathbf{A}_i$, $\psi \opr{|}
\mathbf{A}_i$ are inner conjugate then $\phi$, $\psi$ are inner
conjugate.
\end{prop}
\begin{proof}
For $i \in \mathbb{N}$, let $u_i \in \Sharpo{(\mathfrak{B})}$
intertwine $\phi \opr{|} \mathbf{A}_i$, $\psi \opr{|} \mathbf{A}_i$.
By the lifting theorem for partial isometries, there is a unitary $w_i
\in \mathfrak{B}$ such that $\pi(w_i) = u_i$.  By saturation the
sequence $\{w_i\}_{i \in \mathbb{N}}$ extends to a hyperfinite
sequence $\{w_i\}_{1 \leq i \in N}$ of unitaries in $\mathfrak{B}$.
By applying Robinson's argument (see the proof of
Proposition~\ref{lifting-sequences-of-moorphisms}) and the fact
$\mathbf{A}$ is separable, we conclude there is an $M \cong \infty$
such that $\pi(w_M)^\ast \psi(x) \pi(w_M) = \phi(x)$ for all $x \in
\mathbf{A}$.
\end{proof}

If $\{\mathbf{A}_i\}_{i \in \mathbb{N}}$ is a nondecreasing sequence
of finite dimensional $C^\ast$-algebras with inclusion mappings
$h_k:\mathbf{A}_k \rightarrow \mathbf{A}_{k+1}$, then the sequence of
dimension vectors and inclusion matrices satisfy
$$ \Lambda(h_i) \ \overline{\opr{dim}}(\mathbf{A}_i)=
\overline{\opr{dim}}(\mathbf{A}_{i+1}).$$  
Conversely, any sequence of positive integer vectors $\{\bar{n}_k\}_{k
\in \mathbb{N}}$ and inclusion matrices $\{\Lambda_k\}_{k \in
\mathbb{N}}$ such that $ \Lambda_k \ \bar{n}_k = \bar{n}_{k+1}$ is the
sequence of dimension vectors and inclusion matrices for some
nondecreasing sequence of finite dimensional $C^\ast$-algebras.  The
equivalence of this scheme for describing inductive systems of finite
dimensional $C^\ast$-algebras and Bratteli diagrams~\cite{bratteli} is
well-known.
In what follows, we consider AF algebras with a particular
representation $\mathbf{A} = \overline{\bigcup_{i=1}^\infty
\mathbf{A}_i}$ with $\mathbf{A}_i$ finite dimensional.

If $\mathfrak{B}$ hyperfinite dimensional, the
{\em matrix sequence} of any $C^\ast$-algebra morphism
$\phi:\mathbf{A} \rightarrow \Sharpo{(\mathfrak{B})}$ is the sequence
$\{\Lambda(\Phi \opr{|} \transfer{\mathbf{A}_k})\}_{k \in
\mathbb{N}}$ where $\Phi$ is a lifting of $\phi$ as is guaranteed by
Proposition~\ref{lifting-sequences-of-moorphisms}.  The mapping
matrices are independent of $\Phi$.  In fact,
\begin{equation}\label{infinitesimal-invariance-of-matrix}
\Lambda(\Phi \opr{|} \transfer{\mathbf{A}_k}) = \Lambda(\Phi_k)
\end{equation}
where $\Phi_k: \transfer{\mathbf{A}_k} \rightarrow \mathfrak{B}$ is
any lifting to an internal $C^\ast$-algebra morphism, as follows from:

\begin{lem}
Consider the context of Corollary~\ref{lifting-individual-morphisms}.
The mapping matrix of $\Lambda(\Phi)$ is independent of the lifting
$\Phi$ of $\phi$.
\end{lem}
\begin{proof}
Again assume $\phi$ is injective.  From the lifting result for partial
isometries, there is a unitary $w \in \mathfrak{B}$ such that $w^\ast
\alpha^k_{i \, j} w = \beta^k_{i \, j}$ where $\alpha^k_{i j}$,
$\beta^k_{i j}$ are liftings for $\phi(e_{i j}^k)$. It follows that
all liftings of $\phi$ are conjugate.
\end{proof}

We will denote the matrix sequence of $\phi$ by $\Lambda_\Box(\phi)$.

\begin{prop}\label{bijection-proposition}
$\Lambda_\Box$ determines a bijection from the inner conjugacy classes
of $C^\ast$-morphisms $\mathbf{A} \rightarrow \Sharpo{(\mathfrak{B})}$
and sequences of mapping matrices $\{\Gamma_k\}_{k \in \mathbb{N}}$
such that
\begin{equation} \label{1st-condition-mapping-sequence}
\Gamma_k \ \overline{\opr{dim}}\mathbf{A}_k =
\overline{\opr{dim}}\mathfrak{B}
\end{equation}
and 
\begin{equation} 
\label{2nd-condition-mapping-sequence}
\Gamma_{k+1} \ \Lambda(h_k) = \Gamma_{k}.
\end{equation}
\end{prop}
\begin{proof}
Injectivity follows from Proposition~\ref{conjugacy-condition} and the
facts listed in \S\ref{fd-category}.  We prove $\Lambda_\Box$ is
surjective.  Suppose $\{\Gamma_k\}_{k \in \mathbb{N}}$
satisfies~(\ref{1st-condition-mapping-sequence})
and~(\ref{2nd-condition-mapping-sequence}).
By saturation and overspill,  the sequences $\{\Gamma_k\}_{k \in
\mathbb{N}}$, $\{\transfer{\mathbf{A}_k}\}_{k
\in \mathbb{N}}$ and
$\{\transfer{h_k}\}_{k \in \mathbb{N}}$ extend to hyperfinite
sequences $\{\Gamma_k\}_{1 \leq k \leq M}$, $\{\mathfrak{A}_k\}_{1
\leq k \leq M}$ and $\{\Psi_k\}_{1 \leq k \leq M}$ such that
$\Gamma_k$ is a mapping matrix, $\mathfrak{A}_k$ is a hyperfinite
dimensional $C^\ast$-algebra and
$\Psi_k: \mathfrak{A}_k \rightarrow \mathfrak{A}_{k+1}$ is an
injective morphism of $C^\ast$-algebras satisfying
\begin{equation}\label{divisibility-equation} \Gamma_k \
\overline{\opr{dim}} (\mathfrak{A}_k)= \overline{\opr{dim}}
(\mathfrak{B}).  \end{equation}
There is a $C^\ast$-algebra morphism
$\Phi: \mathfrak{A}_M \rightarrow \mathfrak{B}$ such that
$\Gamma_M = \Lambda(\Phi)$.  For $1 \leq k \leq M$ define the
morphism $\Phi_k : \mathfrak{A}_k \rightarrow \mathfrak{B}$ by
$$\Phi_k = \Phi \ \Psi_{M-1} \ \Psi_{M-2}  \ \cdots \  \Psi_{k}.$$
$\Phi_k$ is tailor-made to satisfy the equation $\Phi_{k+1} \circ
\Psi_k = \Phi_{k}$ for all $k$ such that $1 \leq k \leq M-1$.
Thus $\Sharpo{(\Phi_{k+1})} \circ h_k = \Sharpo{(\Phi_{k})}$ for $k
\in \mathbb{N}$.  
$\Sharpo{(\Psi_k)}$ is thus a compatible sequence of
$C^\ast$-algebra morphisms $\mathbf{A}_k = \Sharpo{(\mathfrak{A}_k)}
\rightarrow \Sharpo{(\mathfrak{B})}$.  \end{proof}


Clearly, a $C^\ast$-morphism $\phi :\mathbf{A} \rightarrow
\Sharpo{(\mathfrak{B})}$ is an imbedding iff each term in the mapping
sequence $\Lambda_\Box(\phi)$ is an inclusion matrix.

\begin{dfn}
If $\bar{m}$, $\bar{n}$ are dimension vectors $\bar{m}$ divides
$\bar{n}$ denoted $\bar{m} \opr{|} \bar{n}$ iff there is an inclusion
matrix $\Gamma$ such that $\bar{n} = \Gamma \bar{m}$.
\end{dfn} 
\begin{rmks}
Note that the matrix $\Gamma$ is not unique, for example:
$$(10) = \begin{bmatrix} 2 \ 2 \ 3 \end{bmatrix} \begin{pmatrix} 1 \\
1 \\ 2 \end{pmatrix} = \begin{bmatrix} 3 \ 3 \ 2 \end{bmatrix}
\begin{pmatrix} 1 \\ 1 \\ 2 \end{pmatrix} .$$
A similar example shows
that $\Gamma \bar{m} = \Gamma \bar{m}'$ does not imply $\bar{m} =  \bar{m}'$.
Consequently, a sequence of dimension vectors $\{\bar{n}_i\}_{i \in
\mathbb{N}}$ such that $n_i \opr{|} n_{i+1}$ does not by itself
determine an increasing sequence of finite dimensional
$C^\ast$-algebras.
\end{rmks}

\begin{prop}\label{imbedding-condition}
Suppose $\{\mathbf{A}_i\}_{i \in \mathbb{N}}$ is a nondecreasing
sequence of finite dimensional $C^\ast$-algebras with inclusion
mappings $h_k:\mathbf{A}_k \rightarrow \mathbf{A}_{k+1}$.  A necessary
and sufficient condition there exist an injective $C^\ast$-algebra
morphism from $\overline{\bigcup_{i=1}^\infty\mathbf{A}_i}$ into a
$C^\ast$-algebra $\Sharpo{(\mathfrak{B})}$, where $\mathfrak{B}$ is a
hyperfinite dimensional $C^\ast$-algebra is that for every $k \in
\mathbb{N}$, $\overline{\opr{dim}} (\mathbf{A}_i)$ divide
$\overline{\opr{dim}} (\mathfrak{B})$.
\end{prop}
\begin{proof}
Necessity: Let $\Lambda_\Box(\phi) = \{\Lambda_i\}_{i \in \mathbb{N}}$.  Then
$\Lambda_i \ \overline{\opr{dim}} (\transfer{\mathbf{A}_\ell}) =
\overline{\opr{dim}} (\mathfrak{B})$.
Thus $ \overline{\opr{dim}} (\transfer{\mathbf{A}_\ell}) =
\overline{\opr{dim}} (\mathbf{A}_\ell)$ divides $\overline{\opr{dim}}
(\mathfrak{B})$.

Sufficiency: For $k \in \mathbb{N}$, let $\Gamma_{k}$ be an inclusion
matrix such that
$$ \Gamma_i \ \overline{\opr{dim}}
(\transfer{\mathbf{A}_\ell}) = \overline{\opr{dim}} (\mathfrak{B}).$$
We claim that $\Gamma_{k}$ can be chosen so that the following
compatibility condition holds:
\begin{equation}\label{compatibility-equation}
 \Gamma_{k+1} \ \Lambda(\psi_k)  =  \Gamma_{k} \quad \mbox{for $1 \leq k \leq M$.}
\end{equation} 
Let $\{\Gamma_k\}_{1 \leq k \leq M}$, $\{\mathfrak{A}_k\}_{1 \leq k
\leq M}$ and $\{\Psi_k\}_{1 \leq k \leq M}$ be as in the proof of
Proposition~(\ref{bijection-proposition}). By overspill, there is no
loss of generality in assuming $\Gamma_k$ is an inclusion matrix for
all $k \leq M$. Define a new hyperfinite sequence $\Gamma'_k$ as
follows:
$$\Gamma'_k = \Gamma_M \ \Lambda(\Psi_{M-1})  \ \Lambda(\Psi_{M-2}) \ \cdots 
\ \Lambda(\Psi_{k+1}) \ \Lambda(\Psi_{k}).$$
$\Gamma'_k$ is an inclusion matrix (composition of inclusion matrices
is an inclusion matrix).  Moreover,
Equation~(\ref{divisibility-equation}) continues to hold since
$$
\begin{aligned}
\Gamma'_k \  \overline{\opr{dim}} (\mathfrak{A}_k) & = \Gamma_M \ \Lambda(\Psi_{M-1}) \
\Lambda(\Psi_{M-2}) \ \cdots \ \Lambda(\Psi_{k+1}) \ \Lambda(\Psi_{k})
\ \overline{\opr{dim}} (\mathfrak{A}_k) \\ 
& = \Gamma_M \ \overline{\opr{dim}} (\mathfrak{A}_M) \\
& = \overline{\opr{dim}} (\mathfrak{B}).
\end{aligned}
$$
Now apply surjectivity of the mapping $\Lambda_\Box$.\end{proof}

Note that the previous result has no reference to the inclusion
matrices $\Lambda(h_k)$. If $\mathbf{A}$ is a UHF algebra
(see~\cite{glimm}), then the sequence of dimension vectors is actually
a sequence of positive integers $\{n_k\}$ such that $n_k \opr{|}
n_{k+1}$. In particular, if $\mathfrak{B}$ is a full internal matrix
algebra of size $N$, $\mathbf{A}$ is injectable into
$\Sharpo{(\mathfrak{B})}$ iff all the integers $n_k \opr{|} N$.

\begin{rem} It follows from the above remarks that if $N$ is a
nonstandard prime, then $\Sharpo{(\opr{M}_N(\transfer{\mathbb{C}}))}$
contains no unitally embedded {\em full} matrix algebra.
\end{rem}

\section{Subhyperfiniteness} 

\begin{dfn}
A $C^\ast$-algebra $\mathbf{A}$ is subhyperfinite iff there is an
internal hyperfinite dimensional $\mathfrak{A}$ such that $\mathbf{A}
\subseteq \Sharpo{(\mathfrak{A})}$.
\end{dfn}
\begin{example} Any approximately finite
$C^\ast$-algebra is subhyperfinite. In particular any commutative
$C^\ast$-algebra is subhyperfinite. \end{example}

\begin{example} If $\mathbf{A}$ is separable residually
finite-dimensional, that is $\mathbf{A}$ has a separating family of
finite dimensional representations, then $\mathbf{A}$ is
subhyperfinite. To show this, note that by separability of
$\mathbf{A}$, we can assume there is a countable separating family
$\phi_k: \mathbf{A} \rightarrow \opr{M}_{n_k}(\mathbb{C})$ of
involutive morphisms.  Let $N \cong \infty$ and consider the
hyperfinite sequence of internal involutive morphisms
$\transfer{\phi_k}: \transfer{\mathbf{A}} \rightarrow
\opr{M}_{\transfer{n}_k}(\transfer{\mathbb{C}})$ for $k \leq N$.
Now
$$\mathfrak{M} = \bigoplus_{k=1}^N \opr{M}_{\transfer{n}_k}(\transfer{\mathbb{C}})$$
is hyperfinite dimensional and $\Phi= \bigoplus_{k=1}^N \transfer{\phi_k}$
is an internal involutive morphism $\transfer{\mathbf{A}} \rightarrow
\mathfrak{M}$.
For $k \in \mathbb{N}$, the following diagram is commutative:
\begin{equation}
\begin{CD}
\transfer{\mathbf{A}}  @>{\Phi}>>   \mathfrak{M} \\
@A{\transfer{}}AA     @VV{\pi \opr{\circ} \opr{proj}_k}V \\
\mathbf{A}   @>{\phi_k}>>   \opr{M}_{n_k}(\mathbb{C})
\end{CD}
\end{equation}	
where $\opr{proj}_k$ denotes the projection of $\mathfrak{M}$ onto
$\opr{M}_{\transfer{n}_k}(\transfer{\mathbb{C}})$.  It follows that
$\Sharpo{(\Phi)}:\Sharpo{(\transfer{\mathbf{A}})} \rightarrow
\Sharpo{(\mathfrak{M})}$ restricted to $\mathbf{A}$ is separating.

For a number of examples of residually finite $C^\ast$-algebras
see~\cite{goodearl-menal}. \end{example}

\begin{example} A $C^\ast$-algebra with proper isometry $w$ (that is,
$\opr{initial \ projection}w = 1$, $\opr{final \ projection}w \neq 1$)
is not subhyperfinite.  To see this, suppose $w \in \mathbf{A}
\subseteq \Sharpo{(\mathfrak{A})}$.  By
Proposition~\ref{lifting-partial-isometries-etc}, there is a partial
isometry $\tilde{w} \in \mathfrak{A}$ such that
$\opr{\pi}{(\tilde{w})} = w$. Since $\pi(\opr{initial \ projection}
\tilde{w}) =1_\mathbf{A}$, by the remarks~\ref{no-index-obstructions},
$\opr{initial \ projection} \tilde{w} =1_\mathfrak{A}$. Thus
$\tilde{w}$ is a proper isometry in $\mathfrak{A}$.  It follows
$\mathfrak{A}$ cannot be hyperfinite dimensional.
\end{example}

We will subsequently provide a standard characterization for
subhyperfinitness.  In order to do this we need some preliminary
results on completely positive contractions.

\section{Complete Positivity} 

In this section we prove a lifting theorem for completely positive
maps very similar to the main result of~\cite{choi-effros}.  The
essential difference between the two results is that in this paper we
have a map between objects which are $C^\ast$-algebras in two
different models of a theory.

We note that complete positivity for a internal linear map $\Phi:
\mathfrak{A} \rightarrow \mathfrak{B}$ between internal
$C^\ast$-algebras is an internal condition.  It means that for every
$n \in \transfer{\mathbb{N}}$, $\opr{M_n}(\Phi):
\opr{M_n}(\mathfrak{A}) \rightarrow \opr{M_n}(\mathfrak{B})$ is
positive.  This implies that for all $n \in \mathbb{N}$,
$\opr{M_n}(\Sharpo{(\Phi)}) : \opr{M_n}(\Sharpo{(\mathfrak{A})})
\rightarrow \opr{M_n}(\Sharpo{(\mathfrak{B})})$ is positive, i.e.,
$\Sharpo{(\Phi)}$ is completely positive.

\begin{dfn}
An internal operator $V$ between internal normed spaces is a {\em near
contraction} iff $\|V\| \leq 1 + \text{infinitesimal}$.
\end{dfn}

\begin{lem}
Suppose $V: \mathcal{H} \rightarrow \mathcal{H}'$ is a near
contraction where $\mathcal{H}$, $\mathcal{H}'$ are internal Hilbert
spaces. Then there is an internal contraction $W:\mathcal{H}
\rightarrow \mathcal{H}'$ such that $V \cong W$ in operator norm.
\end{lem}
\begin{proof}
Consider the self adjoint operator $T= V^\ast V$ on
$\mathcal{H}$.  Let $f$ be the standard
function such that $f(t) =1$ for $t < 1$ and $f(t) = t^{-1}$ for $t
\geq 1$.
By properties of the functional calculus and the fact $T$ is a near
contraction, $S = f(T)$ is such that $S \cong 1$ and $S \, T$ is a
contraction. Then $W = V \, \sqrt{S}$ is a contraction since
$$ W^\ast W =  \sqrt{S} \, V^\ast \, V \, \sqrt{S} = S \, T.$$
Moreover, $W - V = V \, (1 - \sqrt{S}) \cong 0$.
\end{proof}

\begin{prop}\label{perturbation-of-completely-positive-near-contractions}
Suppose $\Phi :\mathfrak{A} \rightarrow \mathfrak{B}$ is a completely
positive near contraction. Then there is a completely positive
contraction $\Psi:\mathfrak{A} \rightarrow \mathfrak{B}$ such that
$\Phi \cong \Psi$ in norm.
\end{prop}
\begin{proof}
By the Stinespring factorization theorem and transfer we can assume
there is a nondegenerate representation $\rho$ of $\mathfrak{A}$ and
an operator $V: \mathcal{H}_\rho \rightarrow \mathcal{H}$ such that
$\Phi(x) = V^\ast\rho(x) V$. Now $ \|\Phi \| = \|V^\ast V \|$.  This
is obvious since we are assuming $\mathfrak{A}$ has an identity. It
follows $V$ is a near contraction.  By the lemma there is a
contraction $W$ such that $W \cong V$.  Thus $\Psi(x) = W^\ast\rho(x)
W$ is a contraction and $\Psi \cong \Phi$ in norm.
\end{proof}

A completely positive contraction $\phi:\mathbf{A} \rightarrow
\mathbf{B}$ is {\em matricial} iff it has a factorization
\begin{equation}\label{matricial-factoring-diagram}
\mathbf{A} \stackrel{\psi}{\longrightarrow} \mathbf{C}
\stackrel{\rho}{\longrightarrow} \mathbf{B} 
\end{equation}
where $\mathbf{C}$ is a finite dimensional $C^\ast$-algeba and $\psi$,
$\rho$ are completely positive contractions.

\begin{prop} \label{lifting-matricial}
Suppose $\mathbf{A}$ is a standard $C^\ast$-algebra, $\mathfrak{B}$ an
internal $C^\ast$-algebra, $\phi:\mathbf{A} \rightarrow
\Sharpo{(\mathfrak{B})}$ a completely positive matricial map.  Then
there is an internal completely positive contraction ${\Phi}:
\transfer{\mathbf{A}} \rightarrow \mathfrak{B}$ which makes the
following diagram commutative:
\begin{equation}\label{matricail-lifting-diagram}
\begin{CD}
\transfer{\mathbf{A}}  @>{{\Phi}}>>   \mathfrak{B} \\
@A{\transfer{}}AA     @VV{\pi}V \\
\mathbf{A}   @>{\phi}>>   \Sharpo{(\mathfrak{B})}
\end{CD}
\end{equation}
\end{prop}
\begin{proof}
By assumption, there is a diagram~(\ref{matricial-factoring-diagram})
with $\mathbf{B} = \Sharpo{(\mathfrak{B})}$.  Now we have the
commutative diagram,
\begin{equation}\label{partial-lifting-diagram}
\begin{CD}
\transfer{\mathbf{A}}  @>{\transfer{\psi}}>>   \transfer{\mathbf{C}}  \\
@A{\transfer{}}AA     @AA{\transfer{}}A \\
\mathbf{A}   @>{\psi}>>   \mathbf{C}
\end{CD}
\end{equation}
To complete the proof, we need the following lemma:

\begin{lem}\label{lifting-cp-maps}
In the context of Proposition~\ref{lifting-matricial}, suppose in
addition $\mathbf{A}$ is a standard finite dimensional
$C^\ast$-algebra.  Then there is an internal completely positive
contraction ${\Phi}: \transfer{\mathbf{A}} \rightarrow \mathfrak{B}$
which makes the diagram~(\ref{matricail-lifting-diagram}) commutative.
\end{lem}
\begin{proof}
$\mathbf{A}$ is the (finite) direct sum of full matrix algebras. The
general case can be reduced to the case $\mathbf{A}$ is a full matrix
algebra by considering the restrictions of $\phi$ to the full matrix
components of $\mathbf{A}$.  Now it is well known that for a
$C^\ast$-algebra $\mathbf{B}$,
$\psi \mapsto \{\psi(e_{i \, j})\}_{i,j}$ 
is a $1$-$1$ correspondence between completely positive maps
$\psi:\opr{M}_n(\mathbb{C}) \rightarrow \mathbf{B}$ and positive
elements of $\opr{M}_n(\mathbf{B})$, where $\{e_{i \, j}\}_{i,j}$ is
the canonical system of matrix units for $\opr{M}_n(\mathbb{C})$.
(See for instance~\cite{choi} where this is shown for the case
$\mathbf{B} = \opr{M}_m(\mathbb{C})$. The general case follows
immediately from this case by considering a faithful involutive
representation of $\mathbf{B}$ on a Hilbert space $H$ and compressing
to finite dimensional subspaces). We claim there is a matrix
$\tilde{B}=\{\tilde{b}_{i \, j}\}_{i,j} \in \opr{M}_n(\mathfrak{B})$
such that
$\opr{M}_n(\pi)(\tilde{B}) = \{\phi(e_{i \, j})\}_{i,j}$.  To see
this,
let $\{\phi(e_{i \, j})\} = T^\ast \ T$ with $T \in
\opr{M}_n(\Sharpo{(\mathfrak{B})})$.  By surjectivity of
$\pi_{\opr{M}_n(\mathfrak{B})}$, there is a $\tilde{T} \in
\opr{M}_n(\mathfrak{B})$ such that $\opr{M}_n(\pi)(\tilde{T}) = T$.
Since $\opr{M}_n(\pi)$ is an involutive morphism of rings, letting
$\tilde{B} = \tilde{T}^\ast \ \tilde{T}$ proves the claim.  Now there
is a completely positive map ${\Phi}$ such that ${\Phi}(\transfer{e_{i
\, j}}) = \tilde{b}_{i \, j}$. By $\mathbb{C}$-linearity of everything
involved, the diagram~(\ref{matricail-lifting-diagram}) is
commutative.  $\Phi$ is a near contraction since $\phi$ is a
contraction.  To complete the proof, apply
Proposition~\ref{perturbation-of-completely-positive-near-contractions}.\end{proof}

To complete the proof of the proposition, instantiate $\mathbf{A}$ in
the lemma with $\mathbf{C}$.  \end{proof}

A completely positive contraction $\phi: \mathbf{A} \rightarrow
\mathbf{B}$ is {\em nuclear} iff there is a sequence $\phi_n$ of
matricial completely positive contractions such that $\phi_n
\rightarrow \phi$ in the point norm topology.  
\begin{prop} \label{lifting-nuclear}
Suppose $\mathbf{A}$ is a standard separable $C^\ast$-algebra,
$\mathfrak{B}$ an internal $C^\ast$-algebra, $\phi:\mathbf{A}
\rightarrow \Sharpo{(\mathfrak{B})}$ a completely positive nuclear
contraction.  Then there is an internal completely positive
contraction ${\Phi}: \transfer{\mathbf{A}} \rightarrow \mathfrak{B}$
which makes~(\ref{matricail-lifting-diagram}) commutative.
\end{prop}
\begin{proof}
Since $\mathbf{A}$ is separable, there is a sequence of completely
positive matricial contractions
$\phi_i:\mathbf{A} \rightarrow \Sharpo{(\mathfrak{B})}$ which
converges to $\phi$ in the point-norm topology. By
Proposition~\ref{lifting-matricial}, each $\phi_i$ has a completely
positive lifting $\Phi_i$ in the sense that the diagram
in~(\ref{matricail-lifting-diagram}) commutes.  The remainder of the
proof shows that the sequence $\{\Phi_i\}_{i \in \mathbb{N}}$ can be
extended in such a way that some $\Phi_N$ with $N \cong \infty$ is the
desired lifting.  Let $\{x_\ell\}_{\ell \in \mathbb{N}}$ be norm dense
in $\mathbf{A}$.
Let $y = \transfer{[\{x_\ell\}_{\ell \in \mathbb{N}}]}$ and $z =
\transfer{[\{\phi(x_\ell)\}_{\ell \in \mathbb{N}}]}$.  Note that in
both these case, the transfer operator is applied to an entire
sequence. In particular for $i \in \mathbb{N}$, $y_i = \transfer{x}_i$
and $z_i = \transfer{[\phi(x_i)]}$. The following formula holds:
\begin{equation} \label{asymptotic-lifting}
\forall k \in \mathbb{N} \ \forall i \in \mathbb{N} \ \exists n  \in \mathbb{N} \
\forall j \geq n \quad \|\phi_j(x_i) - \phi(x_i) \| \leq \frac{1}{k}
\end{equation} 
From this immediately follows:
\begin{equation} \label{asymptotic-lifting-1}
\forall k \ \forall i \ \exists n \quad  \|\Phi_n(y_i) - z_i \| \leq \frac{2}{k}
\end{equation}
By saturation and overspill, the sequence $\{\Phi_\ell\}$ has an
internal extension to a sequence of internal completely positive
contractions $\{\Phi_\ell\}_{1 \leq \ell \leq N}$. By overspill, there
is an $ N_0 \cong \infty$ such that~(\ref{asymptotic-lifting-1})
continues to hold for $k, i \leq N_0$.  Instantiating $k$
in~(\ref{asymptotic-lifting-1}) with some value $\cong \infty$,
$$\forall i \ \exists n_i \  \quad \|\Phi_{n_i}(y_i) - z_i
\| \leq \frac{2}{k} \cong 0.$$
In particular, let $N = \max\{n_1, n_2, \ldots , n_M\}$ where $M \cong \infty$. Then
$$\forall i \in \mathbb{N} \ \|\Phi_N(y_i) - z_i \| \cong 0,$$
which is the desired result. \end{proof} 

\begin{rem}
In the above proposition, if $\phi$ is unital $\Phi$ can be taken to be
unital as well. For $a = \Phi(1_{\transfer{\mathbf{A}}}) \cong
1_{\mathfrak{B}}$.  Replace $\Phi$ by $a^{-1/2} \ \Phi( \cdot) \
a^{-1/2}$.
\end{rem}
\section{Standard Characterization of Subhyperfinite
$C^\ast$-algebras} \label{standard-characterization}

If $x \in L(H)$ and $K \subseteq H$, $\opr{compr}_K x$ denotes the
compression of $x$ to $K$.

\begin{prop} \label{characterization-of-subhyperfinite}
A necessary and sufficient condition a separable nuclear
$C^\ast$-algebra $\mathbf{A}$ be subhyperfinite is that for every
$\epsilon > 0$ and $x_1, \ldots , x_n \in \mathbf{A}$, there is a
representation $\phi$ of $\mathbf{A}$ and a finite dimensional
subspace $K \subseteq H_\phi$ such that for all $i \leq n$,
$\|x_i\| - \epsilon \leq \|\opr{compr}_K \phi(x_i)\|$ and $\| [\phi(x_i)
, \opr{proj}_K ] \| \leq \epsilon$.
\end{prop}
\begin{proof}
Nuclearity is used only for necessity.  Sufficiency: There is a
countable involutive $\mathbb{Q}$-algebra $V$ which is norm dense in
$\mathbf{A}$. Let $\{x_i\}_{i \in \mathbb{N}}$ enumerate $V$.  For
$\ell \in \mathbb{N}$, let $\phi_\ell$ be a representation of
$\mathbf{A}$ and $K_\ell \subseteq H_{\phi_\ell}$ finite dimensional
such that for $i \leq \ell$,
$$\mbox{$\|x_i\| - 1/\ell \leq \|\opr{compr}_{K_\ell}\phi_\ell(x_i)\|$ and $\|
[\phi_\ell(x_i) , \opr{proj}_{K_\ell} ] \| \leq 1/\ell$.}$$  
By transfer there are $\transfer{\mathbb{N}}$ sequences $\{z_\ell\}$,
$\{\Phi_\ell\}$ and $\{\mathcal{K}_\ell\}$ such that for $\ell \in
\transfer{\mathbb{N}}$, $\Phi_\ell$ is an internal $C^\ast$-morphism
of $\transfer{\mathbf{A}}$ and $\mathcal{K}_\ell \subseteq
\mathcal{H}_{\Phi_\ell}$ is a hyperfinite dimensional space such that
for all $\ell \in\transfer{\mathbb{N}}$ and all $i \leq \ell$.
$$\mbox{$\|z_i\| - 1/\ell \leq
\|\opr{compr}_{\mathcal{K}_\ell}\Phi_\ell(z_i)\|$ and $\|
[\Phi_\ell(z_i) , \opr{proj}_{\mathcal{K}_\ell} ] \| \leq 1/\ell$.}$$
Let $\ell \cong \infty$.  Since $\Phi_\ell$ is a representation of
$\transfer{\mathbf{A}}$, $\opr{compr}_{\mathcal{K}_\ell} \Phi_\ell :
\transfer{\mathbf{A}} \rightarrow \mathfrak{L}(\mathcal{K}_\ell)$ is
an internal completely positive contraction.  We claim $\psi:x \mapsto
\pi (\opr{compr}_{\mathcal{K}_\ell} \Phi_\ell(\transfer{x}))$ is an
injective $C^\ast$-morphism $\mathbf{A}\rightarrow \Sharpo{(
\mathfrak{L}(\mathcal{K}_\ell))}$.  It is clear $\psi$ is a completely
positive contraction.
Now for every $x \in V$, $\|[\Phi_\ell(\transfer{x}) ,
\opr{proj}_{\mathcal{K}_\ell} ] \| \cong 0$.
Thus for all $x, y \in V$,
$$\opr{compr}_{\mathcal{K}_\ell} \Phi_\ell(\transfer{(x \ y)}) \cong
\opr{compr}_{\mathcal{K}_\ell} \Phi_\ell(\transfer{x}) \
\opr{compr}_{\mathcal{K}_\ell} \Phi_\ell(\transfer{y}).$$
It follows
$\psi$ is an involutive morphism on $V$ and for every $x \in V$
$$\|x\| = \| \transfer{x}\| \cong
\|\opr{compr}_{\mathcal{K}_\ell}\Phi_\ell(\transfer{x})\| \cong \|\psi(x)\|.$$
By continuity it follows $\psi$ is an isometric $C^\ast$-morphism.

Necessity: Let $\phi: \mathbf{A} \rightarrow \Sharpo{(\mathfrak{B})}$
be an injective $C^\ast$-morphism where $\mathfrak{B}$ is hyperfinite
dimensional and $\Phi:\transfer{\mathbf{A}} \rightarrow \mathfrak{B}$
a unital completely positive lifting in the sense of
Proposition~\ref{lifting-nuclear}; the lifting exists since the
algebra $\mathbf{A}$ and therefore the map $\phi$ are nuclear.
In particular, for all $x, y \in \mathbf{A}$, $\Phi(\transfer{(x \ y)})
\cong \Phi(\transfer{x}) \ \Phi(\transfer{y})$ in the norm of
$\mathfrak{B}$. There is a hyperfinite dimensional $\mathcal{K}$ and
an imbedding $\mathfrak{B} \subseteq \mathfrak{L}(\mathcal{K})$. By
the $\transfer$-version of the Stinespring factorization theorem,
there is a representation $\Psi$ of $\transfer{\mathbf{A}}$ on an
internal Hilbert space $\mathcal{H} \supseteq \mathcal{K}$ such that
$\Phi = \opr{compr}_{\mathcal{K}} \Psi$.
Thus 
$\rho:x \mapsto \opr{compr}_{\mathcal{K}} \Psi (\transfer{x})$ is an
involutive $\mathbb{C}$-morphism $\mathbf{A} \rightarrow
\mathfrak{L}(\mathcal{K})$ modulo $\cong$. This means
$\rho(x \ y) \cong \rho(x) \ \rho(y)$ in the norm of
$\mathfrak{L}(\mathcal{K})$ for all $x, y \in \mathbf{A}$.
From this follows that 
$[\opr{proj}_{\mathcal{K}}, \Psi (\transfer{x})] \cong 0$ for all $x
\in \mathbf{A}$. Proof: If $T$ is an internal selfadjoint operator
such that $\opr{compr}_{\mathcal{K}} T^2 \cong
{\bigl(\opr{compr}_{\mathcal{K}} T \bigr)}^2$ then 
$$
\opr{proj}_{\mathcal{K}} T^2 \opr{proj}_{\mathcal{K}}  =
\opr{proj}_{\mathcal{K}} T (1 - \opr{proj}_{\mathcal{K}}) T
\opr{proj}_{\mathcal{K}} + \opr{proj}_{\mathcal{K}} T
\opr{proj}_{\mathcal{K}} T \opr{proj}_{\mathcal{K}} 
$$
from which follows $\opr{proj}_{\mathcal{K}} T (1 -
\opr{proj}_{\mathcal{K}}) \cong 0$ and thus
$[T, \opr{proj}_{\mathcal{K}}] \cong 0$.
Moreover, 
$$ \|x\| = \|\transfer{x}\| \geq \|\Psi(\transfer{x})\| \geq
\|\opr{compr}_{\mathcal{K}} \Psi(\transfer{x}) \| =
\|\Phi(\transfer{x}) \| \cong \|x\| \quad \mbox{for $x \in
\mathbf{A}$.}$$
Therefore the following formula holds:
$$
\begin{aligned}
\forall \,  &\mbox{standard} \, \epsilon > 0 \ \forall
\, \mathrm{finite} \,F \subseteq \mathbf{A} \ \exists \,
\mbox{internal} \, \Psi \ \exists \, \mbox{internal} \, \mathcal{K} \
\forall z \in \transfer{F} \\
& \begin{array}{cl}
\mathrm{(a)} & \mbox{$\Psi$ is a $C^\ast$-representation of
$\transfer{\mathbf{A}}$,} \\
\mathrm{(b)}  & \mbox{$\mathcal{K} \subseteq \mathcal{H}_\Psi$ is a hyperfinite
dimensional subspace,} \\
\mathrm{(c)} & \mbox{$ \|z\| - \epsilon \leq \|\opr{compr}_{\mathcal{K}} \Psi(z) \|$
and
$\|[\opr{proj}_{\mathcal{K}}, \Psi (z)]\| \leq \epsilon$.}
\end{array}
\end{aligned}
$$
Note that the stronger formula in which the existential quantifiers
are outermost is valid, but is unsuitable for transfer.  By transfer
it follows the same formula is valid with the quantifiers on the first
line ranging over standard values.  $$ \begin{aligned}
\forall & \, \epsilon > 0 \ \forall \, \mathrm{finite} \,F \subseteq
\mathbf{A} \ \exists \,
 \psi \ \exists \, K \
\forall z \in F \\
& \begin{array}{cl}
\mathrm{(a)} & \mbox{$\psi$ is a $C^\ast$-representation of
${\mathbf{A}}$,} \\
\mathrm{(b)}  & \mbox{$K \subseteq H_\psi$ is a finite
dimensional subspace,} \\
\mathrm{(c)} & \mbox{$ \|z\| - \epsilon \leq \|\opr{compr}_{K} \psi(z) \|$
and
$\|[\opr{proj}_{K}, \psi (z)]\| \leq \epsilon$.}
\end{array}
\end{aligned}
$$
This is exactly the condition in the statement of the Proposition.
\end{proof}

\begin{thm}\label{characteriztion-of-qd}
A necessary and sufficient condition a separable nuclear
$C^\ast$-algebra $\mathbf{A}$ be quasidiagonal is it that be
subhyperfinite.
\end{thm}
\begin{proof}
The condition for subhyperfiniteness given in
Proposition~\ref{characterization-of-subhyperfinite} is exactly
condition (iii) in Theorem 1 of ~\cite{voiculescu}.
\end{proof}

\appendix

\section{Ultraproducts of $C^\ast$-algebras}

In certain cases the $C^\ast$-algebras $\Sharpo{(\mathfrak{A})}$ are
exactly the ultraproducts of finite dimensional
$C^\ast$-algebras~(\cite{dacunha-castelle},\S1).  This is remarked in
the paper~\cite{hensvino}, p 19.  For completeness, we sketch the
argument without the machinery of that paper.
We will refer to the map $\transfer: \MDLM \longrightarrow V(\transfer
\mathbb{R})$ constructed in~\S1.2 of~\cite{albeverio-etal} as a {\em
bounded ultrapower embedding}. We recall the context
of~\cite{albeverio-etal}.

Let $\mathcal{U}$ be a free ultrafilter on $\mathbb{N}$. If $\bar{x} =
\{x_i\}_{i \in \mathbb{N}}, \bar{y} = \{y_i\}_{i \in \mathbb{N}}$, define
$$\bar{x} \sim_\mathcal{U} \bar{y} \mbox{ iff }\{i \in \mathbb{N}:
{x}_i ={y}_i \} \in \mathcal{U}.$$ 
$\langle \bar{x} \rangle$ denotes the $\sim_\mathcal{U}$ equivalence
class of $\bar{x}$.  A sequence $\bar{x}$ is {\em bounded} iff for
some $n \in \mathbb{N}$ and all $i \in \mathbb{N}$, ${x}_i \in
V_n(\mathbb{R})$.  The bounded ultrapower of $\MDLM$, denoted
$V^\mathbb{N}(\mathbb{R})/\mathcal{U}$ is
$\{\langle \bar{x} \rangle: \bar{x} \in \MDLM^\mathbb{N} \mbox{
bounded}\}$.
$i:\MDLM \rightarrow V^\mathbb{N}(\mathbb{R})/\mathcal{U}$ is the
mapping
$x \mapsto \langle \mbox{sequence with constant value $x$} \rangle$.
Define a relation $\in_\mathcal{U}$ on
$V^\mathbb{N}(\mathbb{R})/\mathcal{U}$ by
$$\langle\bar{x}\rangle \in_\mathcal{U} \langle\bar{y}\rangle \mbox{
iff } \{i \in \mathbb{N}: {x}_i \in {y}_i\} \in
\mathcal{U}.$$

Note that any sequence $\bar{r}$ with $r_i \in \mathbb{R}$ for all $i
\in \mathbb{N}$ is bounded as an element of
$V^\mathbb{N}(\mathbb{R})$. The set of hyperreals
$\transfer{\mathbb{R}}$ consists of $\langle \bar{r} \rangle$ for
$\bar{r} \in \mathbb{R}^\mathbb{N}$.  The pointwise operations make
$\transfer{\mathbb{R}}$ an ordered nonarchimedean field.  $\langle
\bar{r} \rangle$ is limited iff there is some positive real $M$ such
that $\{i \in \mathbb{N}:|{r}_i| \leq M\} \in \mathcal{U}$.  Letting
$\lim_\mathcal{U}$ denote the generalized limit on $\mathcal{U}$,
$\langle \bar{r} \rangle$ is infinitesimal iff $\lim_\mathcal{U} r_i =
0$.  $\langle \bar{r} \rangle$ is a hyperinteger iff $\{i \in
\mathbb{N}:{r}_i \in \mathbb{N}\} \in \mathcal{U}$.

The relation $\in_\mathcal{U}$ is not the membership relation on
$V^\mathbb{N}(\mathbb{R})/\mathcal{U}$ although it behaves like one.
In fact we can define an injective mapping $j:
V^\mathbb{N}(\mathbb{R})/\mathcal{U} \rightarrow
V(\transfer \mathbb{R})$ which transports $\in_\mathcal{U}$ to the
membership relation on $V(\transfer \mathbb{R})$.  $j$ is defined
recursively, in such a way that $j \opr{|} \transfer{\mathbb{R}}$, is
the identity.  $j$ is referred to as a {\em Mostowski collapsing}.
Following~\cite{albeverio-etal}, the map $\transfer:\MDLM \rightarrow
V(\transfer \mathbb{R})$ is the composition:
$$\MDLM \stackrel{i}{\longrightarrow} V^\mathbb{N}(\mathbb{R})/\mathcal{U}
\stackrel{j}{\longrightarrow}
V(\transfer \mathbb{R}).$$

We point out that it is possible to transport algebraic and analytic
structures from the universe $\MDLM$ to the universe of internal
sets. To do this we cannot apply the mapping $\transfer$ directly.
For instance, the set of metric spaces which are elements $\MDLM$ is
not in $\MDLM$.

We now recall from~\cite{dacunha-castelle} the definition of {\em
bounded ultraproduct} of a sequence of normed spaces $\{E_i\}_{i \in
\mathbb{N}}$ over a free ultrafilter $\mathcal{U}$ on $\mathbb{N}$.
Let $\Pi_0$ be the space of sequences $\bar{x} \in \prod_{i \in
\mathbb{N}} E_i$ for which $\sup_{i \in \mathbb{N}} \|{x}_i\| <
\infty$.  Let $\|\bar{x}\| = \lim_{\mathcal{U}}
\|{x}_i\|$. $\bar{x} \mapsto \|\bar{x}\|$ is a seminorm on $\Pi_0$.
Finally, let
$\prod_{i \in \mathbb{N}} E_i / \mathcal{U}$ be the normed space
$\Pi_0 / \{\bar{x}\in \Pi_0: \|\bar{x}\| = 0\}$.  If $E_i$ are normed
involutive algebras, $\prod_{i \in \mathbb{N}} E_i / \mathcal{U}$ is
also a normed involutive algebra.

\begin{prop}
Let $\mathcal{U}$ be a free ultrafilter on $\mathbb{N}$.  If
$\transfer: \MDLM \longrightarrow V(\transfer \mathbb{R})$ is the
bounded ultrapower embedding over $\mathcal{U}$, then the
$C^\ast$-algebras $\Sharpo{(\mathfrak{A})}$ for $\mathfrak{A}$
hyperfinite dimensional are exactly the ultraproducts over
$\mathcal{U}$ of finite dimensional $C^\ast$-algebras.
\end{prop}
\begin{proof} We sketch the proof.
Suppose $\mathfrak{A}$ is an internal $C^\ast$-algebra.  There is a
bounded sequence $\bar{\mathbf{A}} = \{\mathbf{A}_i\}_{i \in
\mathbb{N}} \in \MDLM^\mathbb{N}$ such that
$j(\langle\bar{\mathbf{A}}\rangle) = \mathfrak{A}$.
It follows $F= \{i \in \mathbb{N}: \mathbf{A}_i \mbox{ is a
$C^\ast$-algebra}\} \in \mathcal{U}$. Modifying the sequence
$\bar{\mathbf{A}}$ on $\complement{F}$ does not change the value of
$\langle\bar{\mathbf{A}}\rangle$, so we may assume $\mathbf{A}_i$ is a
$C^\ast$-algebra for all $i \in \mathbb{N}$.  For $\bar{x} \in
\prod_{i \in \mathbb{N}} \mathbf{A}_i$, the internal norm of
$j(\langle \bar{x} \rangle)$ is the equivalence class $\langle
\{\|{x}_i\|\}_{i \in \mathbb{N}} \rangle$.  $\opr{Fin}(\mathfrak{A})$
consists of elements of the form $j(\langle \bar{x} \rangle)$ with
$\sup_{i \in \mathbb{N}} \| {x}_i \| < \infty$.
The norm of $j(\langle \bar{x} \rangle)$ is infinitesimal iff
$\lim_\mathcal{U} \|{x}_i\| = 0$.  From this it readily
follows that
$\Sharpo{(\mathfrak{A})}$ is isomorphic to $\prod_{i \in \mathbb{N}}
\mathbf{A}_i / \mathcal{U}$

$\mathfrak{A}$ is hyperfinite dimensional iff $F=\{i \in \mathbb{N}:
\mathbf{A}_i \mbox{ is finite dimensional}\} \in
\mathcal{U}$. Modifying $\langle \mathbf{A} \rangle$ on
$\complement{F}$, we may assume $\mathbf{A}_i$ is finite dimensional
for all $i \in \mathbb{N}$.
\end{proof}

\bibliography{../../../nsa.dir/refs}

\begin{thebibliography}{10}

\bibitem{albeverio-etal}
S.~Albeverio, J.E. Fenstad, R.~H\"{o}egh-Kron, and T.~Lindstr\"{o}m.
\newblock {\em Nonstandard Analysis in Stochastic Analysis and Mathematical
  Physics}.
\newblock Academic Press, New York, 1986.

\bibitem{blackadar-kirchberg}
B.~Blackadar and E.~Kirchberg.
\newblock Generalized inductive limits of of finite dimensional
  {$C^\ast$}-algebras.
\newblock {\em Math. Ann.}, 307:343--380, 1997.

\bibitem{bratteli}
O.~Bratteli.
\newblock Inductive limits of of finite dimensional {$C^\ast$}-algebras.
\newblock {\em Trans. Amer. Math. Soc.}, 171:195--234, 1972.

\bibitem{nate-brown}
N.~P. Brown.
\newblock {On Quasidiagonal C*-algebras}, arXiv:math.OA/0008181.

\bibitem{choi}
M.~D. Choi.
\newblock Completely positive maps on complex matrices.
\newblock {\em Linear Alg. Appl.}, 10:285--290, 1975.

\bibitem{choi-effros}
M.~D. Choi and E.~G. Effros.
\newblock The completely positive lifting problem for {$C^\ast$}-algebras.
\newblock {\em Ann. Math.}, 104:585--609, 1976.

\bibitem{dacunha-castelle}
D.~Dacunha-Castelle and J.~L. Krivine.
\newblock Application des ultraproduits \`{a} l'\'{e}tude des espaces et des
  alg\`{e}bres de banach.
\newblock {\em Studia Mathematica}, 41:315--334, 1972.

\bibitem{davidson}
K.~R. Davidson.
\newblock {\em {$C^\ast$}-algebras by Example}, volume~6 of {\em Fields
  Institute Monographs}.
\newblock American Mathematical Society, Providence, 1996.

\bibitem{glimm}
J.~Glimm.
\newblock On a certain class of operator algebras.
\newblock {\em Trans. Am. Math. Soc.}, 95:318--340, 1960.

\bibitem{goodearl-menal}
K.~R. Goodearl and P.~Menal.
\newblock Free and residually finite-dimensional {$C^\ast$}-algebras.
\newblock {\em J. Functional Anal.}, 90:391--410, 1991.

\bibitem{hensvino}
C.~W. Henson and J.~Iovino.
\newblock {\em Analysis and Logic}, chapter Ultraproducts in Analysis.
\newblock London Mathematical Society Lecture Notes series. Cambridge
  University Press, Cambridge, 2002.

\bibitem{hinokuma-ozawa}
T.~Hinokuma and M.~Ozawa.
\newblock Conversion from nonstandard matrix algebras to standard factors of
  type $\mathrm{II}_1$.
\newblock {\em Ill. Journal Math.}, 37(1):1--13, 1993.

\bibitem{hurd-loeb}
A.~Hurd and P.A. Loeb.
\newblock {\em Introduction to Nonstandard Real Analysis}.
\newblock Academic Press, New York, 1987.

\bibitem{jones-sunder}
V.~Jones and V.~S. Sunder.
\newblock {\em Introduction to Subfactors}, volume 234 of {\em London
  Mathematical Society, Lecture Note Series}.
\newblock Cambridge University Press, 1997.

\bibitem{keisler}
H.J. Keisler.
\newblock {\em An Infinitesimal Approach to Stochastic Analysis}, volume 297 of
  {\em Memoirs, American Mathematical Society}.
\newblock American Mathematical Society, 1984.

\bibitem{nelsonprob}
E.~Nelson.
\newblock {\em Radically Elementary Probability Theory}.
\newblock Princeton University Press, 1987.

\bibitem{rordam-larsen-laustsen}
M.~R{\o}rdam, F.~Larsen, and N.~J. Laustsen.
\newblock {\em An Introduction to {$K$}-Theory for {$C^\ast$}-algebras}.
\newblock London Mathematical Society, Student texts. Cambridge University
  Press, 2000.

\bibitem{jt1}
F.~J. Thayer.
\newblock Nonstandard analysis of graphs.
\newblock {\em Houston Journal of Math.}, To Appear.

\bibitem{voiculescu}
D.~Voiculescu.
\newblock A note on quasidiagonal {$C^\ast$}-algebras and homotopy.
\newblock {\em Duke Math. Journal}, 62(2):267--271, 1991.

\end{thebibliography}
\bibliographystyle{hplain}

\end{document}